\documentclass[12pt]{article}

\newtheorem{df}{Definition}
\newtheorem{thm}[df]{Theorem}
\newtheorem{prop}[df]{Proposition}
\newtheorem{lem}[df]{Lemma}
\newtheorem{rem}[df]{Remark}

\newtheorem{ex}[df]{Example}

\newcommand{\qed}{\hfill$\Box$}

\newcommand{\cir}{{\rm circ}}

\newcommand{\head}{{\frak h}}
\newcommand{\tail}{{\frak t}}

\usepackage{amssymb}
\usepackage[final]{graphicx}
\usepackage[dvipdfmx]{color}
\usepackage{comment}

\title{On the Ihara expression for the generalized \\weighted zeta function}
\author{
	Ayaka Ishikawa\thanks{a-iskw@fc.ritsumei.ac.jp}\\
	Ritsumeikan University
	\\Noji-higashi, Kusatsu 525-85771, Japan	
	\and
	Hideaki Morita\thanks{morita@mmm.muroran-it.ac.jp}\\
	Muroran Institute of Technology\\
	Mizumoto, Muroran 050-8585, Japan
}
\date{}

\begin{document}

%
%

\maketitle

\begin{abstract}
	We consider the generalized weighted zeta function for a finite digraph, 
	and show that it has the Ihara expression, 
	a determinant expression of graph zeta functions, 
	with a certain specified definition for inverse arcs. 
	A finite digraph in this paper allows multi-arcs or multi-loops. 
\end{abstract}

	\section{Introduction}
	
	A graph zeta function is a formal power series associated with a finite graph. 
	It enumerates the closed paths of a given length, exposes the primes, or depicts the cycles 
	in a finite graph. 
	The prototype of the graph zeta functions was introduced by Y. Ihara \cite{ihara66} in 1966 
	from a number-theoretical point of view. 
	Ihara's zeta was subsequently pointed out by J. -P. Serre \cite{serre80} that 
	it can be formulated in terms of finite graphs, 
	and is now called the Ihara zeta function for a finite graph \cite{hashimoto90, kotanisunada00, sunada88}. 
	In the paper of H. Bass \cite{bass92}, 
	as was implied by Ihara \cite{ihara66}, 
	the Ihara zeta is provided the determinant expression 
	described by the adjacency matrix and the degree matrix of the corresponding graph. 
	This determinant expression is now called the Ihara expression 
	\cite{IIMSS21, ishikawamoritasato22+} (see also \cite{morita20}), 
	and the theorem is called the Bass-Ihara theorem. 
	The Ihara expression is one of the main interests in the study of graph zeta functions, 
	and many researches have pursued this subject 
	\cite{bartholdi99, choekwakparksato07, foatazeilberger99, IIMSS21, ishikawamoritasato22+, mizunosato04, 
	northshield99, sato07, starkterras96} . 
	It is also necessary to mention that 
	the Ihara expression was recently provided a new point of view from quantum walk theory 
	\cite{ishikawa22+, konnosato12}, 
	and thus the significance of the Ihara expression is now increasing in related areas.

	The subject of the present paper is the Ihara expression for the generalized weighted zeta function. 
	The generalized weighted zeta function was introduced in \cite{morita20} 
	as a single scheme which unifies the graph zetas appeared in previous studies, 
	for instance, the Ihara zeta \cite{ihara66}, the Bartholdi zeta \cite{bartholdi99}, 
	the Mizuno-Sato zeta \cite{mizunosato04} and the Sato zeta \cite{sato07}. 	
	Graph zetas may have in general four expressions called 
	the exponential expression, the Euler expression, the Hashimoto expression and the Ihara expression. 
	It is verified in \cite{morita20} that 
	the first three expressions are equivalent for those graph zeta functions. 
	The last two expressions are both determinant expressions, 
	where the size of the matrices used in the latter one, the adjacency matrix and the degree matrix, 
	is smaller in general than the other one, the edge matrix.

	For known examples of graph zetas, 
	the Ihara expression is obtained by transforming the Hashimoto expression. 
	In this paper, 
	we will show that 
	the generalized weighted zeta function also have the Ihara expression 
	as in the same manner with those graph zetas. 
	In particular, we verify the main theorem for the case where 
	the underlying graph is an arbitrary finite digraph. 
	Graph zetas have usually been defined via the symmetric digraph of a given finite graph, 
	so it is natural to define for finite digraphs rather than finite graphs. 
	In addition, a finite digraph in this paper allows multi-arcs and multi-loops, 
	and one will see in the procedure that it is an unavoidable issue 
	how one defines the inverses for each arc of a digraph. 
	For this, we can consider two extreme ways. 
	One is the case where all the arcs with inverse direction to an arc $a$ are defined to be the inverses of $a$, 
	and the other one is the case where a single arc with inverse direction, if exists, defined to be the inverse of $a$. 
	In \cite{ishikawamoritasato22+}, we treat the former case. 
	In the present article, we treat the latter case. 
	These two cases are natural generalizations for the case where 
	the underlying graph is a finite simple graph. 
	Therefore, 
	the main result in this paper provides a way to generalize the developments in previous researches on the Ihara expression, 
	and gives a unified method to handle it.  
	
	Throughout this paper, 
	we use the following notation. 
	The ring of integers is denoted by ${\mathbb Z}$. 
	The field of rational numbers and complex numbers are denoted by 
	${\mathbb Q}$ and ${\mathbb C}$ respectively. 
	For a set $X$, the cardinality of $X$ is denoted by $|X|$. 
	The Kronecker delta is denoted by $\delta_{xy}$, 
	which returns $1$ if $x=y$, $0$ otherwise. 
	The symbol $I$ stands for the identity matrix.

	\section{Preliminaries}\label{section : Preliminaries} 
	
		\subsection{Graphs and Digraphs}\label{subsection : Graphs and Digraphs}
	
		A {\it digraph} is a pair $\Delta=(V, {\cal A})$ of a set $V$ and 
		a multi-set ${\cal A}$ consisting of ordered pairs $(u, v)$ of elements $u ,v$ in $V$. 
		If the cardinalities of $V$ and ${\cal A}$ are finite, then $\Delta$ is called a {\it finite} digraph. 
		An element of $V$ (resp. ${\cal A}$) is called a {\it vertex} (resp. an {\it arc}) of  $\Delta$. 
		An arc $a=(u, v)$ is depicted by an arrow from $u$ to $v$. 
		The vertex $u$ is called the {\it tail} of $a$, 
		and $v$ the {\it head} of $a$, 
		which are denoted by $\tail(a)$ and $\head(a)$ respectively. 
		Let $u, v\in V$. 
		We denote by ${\cal A}_{uv}$ the set 
		$$
			\{a\in{\cal A}\mid \tail(a)=u, \head(a)=v\}
		$$
		of arcs with the tail $u$ and the head $v$. 
		An arc $l$ of the form $(u, u)$ is called a {\it loop}. 
		The set of loops is denoted by ${\cal L}$. 
		Hence ${\cal L}=\sqcup_{u\in V}{\cal A}_{uu}$. 
		Note that $|{\cal A}_{uv}|\geq 1$ may occur in general. 
		Thus an arc $a\in{\cal A}_{uv}$ sometimes called a {\it multi-arc}, 
		and the cardinality $|{\cal A}_{uv}|$ is called the {\it multiplicity} of $a$. 
		Similarly, a loop $l\in {\cal A}_{uu}$ may be called a {\it multi-loop}. 
		The cardinality $|A_{uu}|$ is called the {\it multiplicity} of $l$. 
		Set 
		$
			{\cal A}_{u*}=\sqcup_{v\in V}{\cal A}_{uv}
		$ 
		and 
		$
			{\cal A}_{*v}=\sqcup_{u\in V}{\cal A}_{uv}.
		$ 
		A digraph $\Delta=(V, {\cal A})$ is called {\it simple} 
		if ${\cal A}_{uu}=\emptyset$ for any $u\in V$ and $|{\cal A}_{uv}|=1$ if ${\cal A}_{uv}\neq\emptyset$. 
		Let ${\cal A}(u, v)={\cal A}_{uv}\cup{\cal A}_{vu}$ denote 
		the set of arcs lying between vertices $u$ and $v$. 
		A digraph $\Delta$ is called {\it connected} 
		if ${\cal A}(u ,v)\neq\emptyset$ for any distinct $u$, $v$. 	
		A digraph in this paper is always assumed to be connect otherwise stated.

		Let $\Delta=(V, {\cal A})$ be a finite digraph, and $u$, $v$ two distinct vertices. 
		We may assume that $|{\cal A}_{uv}|\leq |{\cal A}_{vu}|$. 
		If ${\cal A}_{uv}$ and ${\cal A}_{vu}$ are both not empty, 
		then one can fix an injection 
		$$\iota_{uv} : {\cal A}_{uv}\rightarrow {\cal A}_{vu}.$$ 
		In this case, we say that an arc $a\in {\cal A}_{uv}$ {\it has inverse}, 
		and the arc $\iota_{uv}(a)\in {\cal A}_{vu}$ is the {\it inverse arc}, or simple the {\it inverse} of $a$, 
		denoted by $a^{-1}$; 
		and vice versa, $a$ is the inverse of $a^{-1}$. 
		An arc $a'\in{\cal A}_{vu}$ not lying in the image of $\iota_{uv}$ has no inverse. 
		In the case where ${\cal A}_{uv}=\emptyset$, 
		any arc $a'\in {\cal A}_{vu}$ is defined to have no inverse. 
		If $u=v$, then $\iota_{uu}$ is defined to be the identity map. 
		By this definition, 
		a loop $l\in{\cal A}_{uu}$ satisfies $l^{-1}=l$, 
		that is, each loop is {\it self-inverse}. 
		Alternatively, one can also define any arc belonging to ${\cal A}_{vu}$ to be inverse of an arc of ${\cal A}_{uv}$. 
		This alternate definition also works, 
		and the development with this definition will be found in \cite{ishikawamoritasato22+}.

		A {\it graph} is a pair $\Gamma=(V, E)$ of a set $V$ and a multi-set $E$ consisting of 
		$2$-subsets $\{u, v\}$ of $V$. 
		If $V$ and $E$ are finite (multi-)sets, 
		then the graph $\Gamma$ is called {\it finite}. 
		An element $\{u, v\}\in E$ is called a {\it edge}. 
		In particular, an edge of the form $l=\{u, u\}$ is called an {\it loop}. 
		The set of loops is denoted by $L$. 
		Obviously, $\{u, v\}\in E\setminus L$ implies $u\neq v$. 
		We also suppose that an edge or a loop has multiplicity. 
		Hence these are sometimes called 
		an {\it multi-edge} and {\it multi-loop}, respectively. 
		In other words, 
		if we denote by $E(u, v)$ the set of multi-edges lying between 
		vertices $u, v\in V$, 
		then we assume that $|E(u, v)|\geq 1$ for $u, v\in V$ with $E(u, v)\neq\emptyset$. 
		Note that $E(u, u)$ denotes the set of loops of the form $\{u, u\}$. 
		The cardinality $|E(u, v)|$ is called the {\it multiplicity} of an edge $\{u, v\}$. 
		A graph is called {\it simple} if 
		it has no loops and the multiplicity of any edge is at most one. 
		The matrix 
		$$A_{\Gamma}=(|E(u, v)|)_{u, v\in V}$$ 
		is called the {\it adjacency matrix} of $\Gamma$. 
		For a vertex $u\in V$, the number of edges $\{u, v\}$ ($v\in V$) is called 
		the {\it degree} of $u$, denoted by $d_u$. 
		Thus we have 
		$
			d_u
			=
			\sum_{v\in V}|E(u, v)|
		$ 
		for $u\in V$. 
		The diagonal matrix 
		$$D_{\Gamma}=(\delta_{uv}d_u)_{u, v\in V}$$ 
		is called the {\it degree matrix} of $\Gamma$.

		Let $\Gamma=(V, E)$ be a finite graph. 
		We recall the definition of the symmetric digraph of $\Gamma$. 
		We assign for each edge $\{u, v\}\in E\setminus L$, 
		two arcs $(u, v)$ and $(v, u)$ in mutually reverse direction. 
		For a loop $\{u, u\}\in L$, we assign a single directed loop $(u, u)$. 
		Then we have an set of arcs 
		$${\cal A}=\{(u, v), (v, u)\mid \{u, v\}\in E\setminus L\}\sqcup\{(u, u)\mid \{u, u\}\in L\}.$$
		The finite digraph constructed in this manner is called the {\it symmetric digraph} 
		of a finite graph $\Gamma$. 
		An arc $a'=(v, u)$ is called the {\it inverse} of $a=(u, v)$, 
		and denote it by $a'=a^{-1}$. 
		Any loop $l=(u, u)$ is defined to be self-inverse, i.e., $l^{-1}=l$. 
		Therefore, the notion of inverse arcs is straightforwardly defined for the symmetric digraph of a finite graph, 
		and one can see that the preceding definitions of inverse arcs, 
		including the alternating one considered in \cite{ishikawamoritasato22+}, are natural generalization 
		of the case for the symmetric digraph. 
		It can readily be confirmed that 
		the symmetric digraph of a simple graph is simple.

		\subsection{The generalized weighted zeta function}\label{subsection : The generalized weighted zeta function}
		
		In this subsection, 
		we briefly review the definition of the generalized weighted zeta function following \cite{morita20}. 
		Let $\Delta=(V, {\cal A})$ be a finite digraph, 
		${\cal A}^{\mathbb Z}$ the set of two-sided infinite sequence. 
		Let $\varphi$ be the left shift operator on ${\cal A}^{\mathbb Z}$, 
		and 
		$$
			\Pi_{\Delta}
			=
			\{
			(a_i)_{i\in{\mathbb Z}}\in{\cal A}^{\mathbb Z}
			\mid
			\head(a_i)=\tail(a_{i+1}),\ \forall i\in{\mathbb Z}
			\}
		$$ 
		the subshift of the dynamical system $({\cal A}^{\mathbb Z}, \varphi)$ 
		consisiting of two-sided infinite path of $\Delta$. 
		If we denote the restriction $\varphi|_{\Xi}$ by $\lambda$, 
		then we have a {\it quasi-finite} dynamical system $(\Pi_{\Delta}, \lambda)$, 
		that is, the set 
		$
			X_m
			=
			\{
			x\in\Pi_{\Delta}
			\mid
			\lambda^m(x)=x
			\}
		$ 
		of $m$-periodic points in $(\Pi_{\Delta}, \lambda)$ 
		is a finite set for each $m\geq 1$, 
		since we have $|X_m|\leq |{\cal A}|^m$. 
		For $x\in X_m$, the integer $m$ is called a {\it period} of $x$. 
		Thus the union $X=\cup_{m\geq 1}X_m$ consists of all periodic points in $(\Pi_{\Delta}, \lambda)$. 
		An element $x\in X_m$ is called a {\it closed path} of $\Delta$ of {\it length} $m$. 
		Note that the union is not disjoint, 
		since any multiple of a period of $x\in X$ is again a period of it. 
		Let $x=(a_i)\in X$, 
		and let a period $m$ of $x$ be fixed. 
		A consecutive $m$-section $(a_k, a_{k+1}, \dots, a_{k+m-1})$ 
		is called a {\it fundamental section} of $x\in X_m$. 
		Thus the fundamental sections of $X_m$ are in one-to-one correspondence 
		with the closed paths of $\Delta$. 
		If we consider the $\varphi$-stable subset 
		$
			\Pi_{\Delta}^{\flat}
			=
			\{
			(a_i)\in\Pi_{\Delta}
			\mid 
			a_i^{-1}\neq a_{i+1},\ \forall i
			\}
		$, 
		then an element of $X_m$ is called 
		a {\it reduced} closed path of length $m$. 
		Let $\varpi(x)$ denote the minimum period of $x$. 
		Hence it is obvious that $x\in X_{\varpi(x)}$ for any $x\in X$. 
		If we consider $x\in X$ as an element of $X_{\varpi(x)}$, 
		then $x$ is called a {\it prime element}. 
		If $x\in X$ is prime, then we denote it by $\pi(x)$. 
		Obviously we have $\pi(x)\in X_{\varpi(x)}$. 
		For the subshift $\Pi_{\Delta}$ (resp. $\Pi_{\Delta}^{\flat}$), 
		a prime element is called 
		a {\it prime} (resp. {\it prime reduced}) closed path of $\Delta$.

		Let $\Delta=(V, {\cal A})$ be a finite digraph, $R$ a commutative ${\mathbb Q}$-algebra. 
		For two functions 
		$\tau, \theta : {\cal A}\rightarrow R$, 
		we consider the following function with two variables 
		$$
		\theta^{\rm G} : {\cal A}\times{\cal A}\rightarrow R
		$$ defined by 
		$
		\theta^{\rm G}(a, a')
		=
		\tau(a')\delta_{\head(a)\tail(a')}
		-\upsilon(a')\delta_{a^{-1}a'}.
		$ 
		If $a$ has no inverse, we understand that $\delta_{a^{-1}a'}$ is zero. 
		Let $x=(a_i)\in\Pi_{\Delta}$ be periodic of period $m$, 
		and $(a_k, a_{k+1}, \dots, a_{k+m-1})$ a fundamental section. 
		The following product 
		$$
			\theta^{\rm G}(\alpha_k, \alpha_{k+1})\theta^{\rm G}(\alpha_{k+1}, \alpha_{k+2})
			\cdots
			\theta^{\rm G}(\alpha_{k+m-2}, \alpha_k+{m-1})\theta^{\rm G}(\alpha_{k+m-1}, \alpha_k)
		$$ 
		is called the {\it circular product} of $\theta^{\rm G}$ {\it along with} $x$. 
		The circular product does not depend the choice of fundamental sections, 
		and it is denoted by $\cir_{\theta^{\rm G}}(x)$,
		Let $N_m(\theta^{\rm G})$ denote the sum 
		$
			\sum_{x\in X_m}\cir_{\theta^{\rm G}}(x).
		$

		\begin{df}[Generalized weighted zeta function]\label{df : Generalized weighted zeta function}
		{\em 
		Let $t$ be an indeterminate. 
		The following formal power series
		$$
			\exp\left[
				\sum_{m\geq 1}
				\frac{N_m(\theta^{\rm G})}{m}t^m
			\right]\in R[[t]], 
		$$
		denoted by $Z_{\Delta}(t; \theta^{\rm G})$, is called the {\it generalized weighted zeta function} for $\Delta$. 
		}
		\end{df}


		\begin{ex}[Ihara zeta function etc]\label{ex : Ihara zeta function etc}
		{\em 
		If we let $\tau=\upsilon=1$, i.e., 
		given a map $\theta^{\rm I} : {\cal A}\times{\cal A}\rightarrow R$ 
		by 
		$$
			\theta^{\rm I}(a, a')
			=\delta_{\head(a)\tail(a')}-\delta_{a^{-1}a'},
		$$ 
		then the generalized weighted zeta function $Z_{\Delta}(t; \theta^{\rm I})$ 
		is called the {\it Ihara zeta function}\cite{bass92, hashimoto90, ihara66, kotanisunada00, serre80}. 
		We have $N_m(\theta^{\rm I})=|X_m^{\flat}|$, 
		the number of reduced closed paths of length $m$ in $\Delta$. 
		If $\tau=1$ and $\upsilon=0$, 
		then $Z_{\Delta}(t; \theta^{\rm G})$ is called the {\it Bowen-Lanford zeta function}\cite{bowenlanford70}. 
		In this case, $N_m(\theta^{\rm G})$ gives the number $|X_m|$ of closed paths of length $m$ in $\Delta$. 
		The {\it Bartholdi zeta function}\cite{bartholdi99} is the case where $\tau=1$ and $\upsilon=1-q$, 
		where $q$ is indeterminate. 
		Similarly, 
		the {\it Mizuno-Sato zeta function}\cite{mizunosato04} is the case where $\tau=\upsilon$. 
		The {\it Sato zeta function}\cite{sato07} is the case where $\upsilon=0$. 
		}
		\end{ex}

		\begin{rem}
		{\em
		If a digraph $\Delta$ consists of several connected components 
		$\Delta=\sqcup_{i=1}^n\Delta_i$, 
		then one can easily see that 
		$
			Z_{\Delta}(t; \theta^{\rm G})
			=
			\prod_{i=1}^n
			Z_{\Delta_i}(t; \theta^{\rm G}).
		$
		}
		\end{rem}

		\subsection{Three expressions}

		The generalized weighted zeta functions have two other expressions. 
		Let $\Delta=(V, {\cal A})$ be a finite digraph, $R$ a commutative ${\mathbb Q}$-algebra, 
		and $\theta : {\cal A}\times{\cal A}\rightarrow R$ a map. 
		Two elements $x, y\in\Pi_{\Delta}$ is called {\it equivalent} iff 
		there exists an integer $k$ satisfying $y=\lambda^k(x)$. 
		We denote by $x\sim y$ this equivalence relation on $\Pi_{\Delta}$. 
		Note that the relation $\sim$ also affords an equivalence relation on $X$. 
		An equivalence class with representative $x\in \Pi_{\Delta}$ is denoted by $[x]$. 
		An element of the coset ${\frak X}=X/\sim$ is called a {\it cycle} of $\Delta$. 
		Since the relation $\sim$ affords an equivalence relation on each $X_m$, 
		we have ${\frak X}=\cup_{m\geq 1}{\frak X}_m$, 
		where ${\frak X}_m=X_m/\sim.$ 
		If $[x]\in{\frak X}_m$, 
		then positive integer $m$ is called the {\it period} of the cycle $[x]$. 
		A cycle $[x]$ with reduced (resp. prime) $x$ is called 
		a {\it reduced} (resp. {\it prime}) cycle of $\Delta$. 
		If $[x]$ is prime, 
		then we denote it by $\pi([x])$, which belongs to ${\frak X}_{\varpi(x)}$. 
		In other words, we have $\varpi([x])=\varpi(x)$. 
		Let $M_{\Delta}(\theta^{\rm G})=(\theta^{\rm G}(a, a'))_{a, a'\in{\cal A}}$,  
		which is a square matrix of degree $|{\cal A}|$. 
		We consider the following two formal power series belonging to $R[[t]]$: 
		
		$$
		E_{\Delta}(t; \theta^{\rm G})=
		\prod_{[x]\in{\frak X}}
		\frac{1}{1-\cir_{\theta^{\rm G}}(\pi([x]))t^{\varpi([x])}}, 
		\quad
		H_{\Delta}(t; \theta^{\rm G})=
		\frac{1}{\det(I-tM_{\Delta}(\theta^{\rm G}))}.
		$$
		
		In \cite{morita20}, 
		it is verified that 
		these three formal power series are identical for the generalized weighted zeta function. 
		
		\begin{prop}
		For a finite digraph $\Delta$, 
		it follows that
		$
		Z_{\Delta}(t; \theta^{\rm G})=E_{\Delta}(t; \theta^{\rm G})=H_{\Delta}(t; \theta^{\rm G}).
		$
		\end{prop}
		
		In general, 
		these identities are not necessarily hold. 
		In a general setting \cite{morita20}, 
		it only holds that 
		$H$ can be reformulate in the form of $E$ by the Foata-Zeilberger theorem\cite{foatazeilberger99} 
		(see \cite{cartierfoata69, foatahan07, konvalinkapak07} for related topics), 
		and also $E$ to $Z$ by simple calculation. 
		The other implications, $Z$ to $E$ and $E$ to $H$, need certain conditions, 
		called the {\lq Euler condition\rq} and the \lq Hashimoto condition\rq. 
		See \cite{morita20} for precise information. 
		On the other hand, it is also shown in \cite{morita20} that, for the generalized weighted zeta, 
		these three expressions are equivalent to each other. 
		These three expressions for the generalized weighted zeta function are called 
		the {\it exponential expression}, the {\it Euler expression} 
		and the {\it Hashimoto expression}, respectively. 
		In particular, the existence of the Hashimoto expression is significant for our development. 
		We will construct the Ihara expression by reformulating the Hashimoto expression 
		(c.f., \cite{watanabefukumizu10}).

		\section{Main result}
		\label{section : Main result}

		Let $\Delta=(V, {\cal A})$ be a finite digraph which allows multi-arcs and multi-loops, 
		and $R$ a commutative ${\mathbb Q}$-algebra. 
		Given two functions $\tau$, $\upsilon : {\cal A}\rightarrow R$, 
		let $\theta^{\rm G}:{\cal A}\times{\cal A}\rightarrow R$ be the map defined by 
		$$
			\theta^{\rm G}(a, a')=\tau(a')\delta_{\head(a)\tail(a')}-\upsilon(a')\delta_{a^{-1}a'},
		$$
		for $a, a'\in{\cal A}$. 
		We consider the generalized weighted zeta function $Z_{\Delta}(t; \theta^{\rm G})$. 
		We are in position to construct the Ihara expression for the generalized weighted zeta 
		$Z_{\Delta}(t; \theta^{\rm G})$ with the definition of inverse arc given in \ref{subsection : Graphs and Digraphs}.

		For $u, v\in V$, recall that 
		$
			\Phi_{\Delta}=
			\{
			(u, v)\in\Phi_{\Delta}
			\mid
			{\cal A}(u, v)\neq\emptyset
			\}.
		$ 
		Thus it follows that ${\cal A}=\sqcup_{(u, v)\in\Phi_{\Delta}}{\cal A}(u, v)$. 
		For vertices $u, v\in V$, 
		we write 
		$
			u\preceq v
		$ 
		if $u\neq v$ and 
		$
			|{\cal A}_{uv}|\leq|{\cal A}_{vu}|.
		$ 
		If $u\neq v$ and 
		$
			|{\cal A}_{uv}|<|{\cal A}_{vu}|,
		$ 
		we write $u\prec v$. 
		Thus, any $(u, v)\in V\times V$ satisfies 
		$u\preceq v$, $v\prec u$ or $u=v$. 
		If $u=v$ or $u\preceq v$, 
		then we fix a injection 
		$
			\iota_{uv} : {\cal A}_{uv}\rightarrow{\cal A}_{vu}
		$ 
		as in section \ref{subsection : Graphs and Digraphs}. 
		In particular, the map $\iota_{uu}$ is assume to be the identity map on ${\cal A}_{uu}$. 
		Let 
		\begin{eqnarray*}
		&&
			\Phi_{\Delta}^{(1)}
			=
			\{
			(u, v)\in\Phi_{\Delta}
			\mid
			u\preceq v, {\cal A}_{uv}\neq\emptyset
			\},\\
		&&
			\Phi_{\Delta}^{(2)}
			=
			\{
			(u, u)\in\Phi_{\Delta}
			\mid
			{\cal A}_{uu}\neq\emptyset
			\},\\
		&&
			\Phi_{\Delta}^{(3)}
			=
			\{
			(u, v)\in\Phi_{\Delta}
			\mid
			u\prec v, {\cal A}_{uv}=\emptyset
			\},
		\end{eqnarray*}
		and let 
		\begin{eqnarray*}
			{\cal A}^{(1)}
			=
			\bigcup_{(u, v)\in\Phi_{\Delta}^{(1)}}
			{\cal A}_{uv},
		& 
			{\displaystyle
			{\cal A}^{(-1)}
			=
			\bigcup_{(u, v)\in\Phi_{\Delta}^{(1)}}
			{\cal A}_{uv}^{-1},
			}
		& 
			\overline{{\cal A}^{(1)}}
			=
			\bigcup_{(u, v)\in\Phi_{\Delta}^{(1)}}
			{\cal A}_{vu}\setminus{\cal A}_{uv}^{-1},
		\\ 
			{\cal A}^{(2)}
			=
			\bigcup_{(u, u)\in\Phi_{\Delta}^{(2)}}
			{\cal A}_{uu},
		& 
			{\displaystyle
			{\cal A}^{(3)}
			=
			\bigcup_{(u, v)\in\Phi_{\Delta}^{(3)}}
			{\cal A}_{vu},
			}
		&
		\end{eqnarray*}
		where 
		$
			{\cal A}_{uv}^{-1}
			=
			\{
			\iota_{uv}(a)\mid a\in{\cal A}_{uv}
			\}.
		$
		Thus, it follows that 
		$
			{\cal A}
			=
			{\cal A}^{(1)}\sqcup{\cal A}^{(-1)}\sqcup\overline{{\cal A}^{(1)}}\sqcup{\cal A}^{(2)}\sqcup{\cal A}^{(3)}. 
		$ 
		The set of arcs with inverse is given by ${\cal A}^{(1)}\sqcup{\cal A}^{(2)}$, 
		and $\overline{{\cal A}^{(1)}}\sqcup{\cal A}^{(3)}$ give the set of arcs without inverse. 
		We set ${\cal A}^{\times}=\overline{{\cal A}^{(1)}}\sqcup{\cal A}^{(3)}$. 
		For an arc $a\in{\cal A}$, 
		let 
		$$
			{\cal E}(a)
			=
			\left\{
				\begin{array}{ll}
				\{a, a^{-1}\},& \mbox{if $a\in{\cal A}^{(1)}$},\\
				\{a\}, & \mbox{otherwise.}
				\end{array}
			\right.
		$$
		and let $c_a(t)=c_a(t; \theta^{\rm G})=1-\prod_{\alpha\in{\cal E}(a)}(-\upsilon(\alpha)t)$. 
		For $u, v\in V$, 
		define
		$$
			a_{uv}
			=
			\sum_{a\in{\cal A}_{uv}}
			\frac{\tau(a)}{c_a(t)}
			\in R[[t]],
			\quad
			b_{uv}
			=
			\delta_{uv}
			\sum_{a\in{\cal A}^{(1)}\cap{\cal A}_{u*}}
			\frac{\tau(a)\upsilon(a^{-1})}{c_a(t)}
			\in R[[t]].
		$$

		\begin{df}
		{\em 
		Let $\Delta=(V, {\cal A})$ be a finite digraph. 
		The following $|V|\times |V|$ matrices 
		$$
			A_{\Delta}(\theta^{\rm G})=(a_{uv})_{u, v\in V}, \quad 
			B_{\Delta}(\theta^{\rm G})=(b_{uv})_{u, v\in V}
		$$ are called the {\it weighted adjacency matrix} and the {\it weighted backtrack matrix} 
		for $\Delta$ respectively.
		}
		\end{df}

		\begin{ex}\label{example : Example}
		{\em
			Let $V=\{1, 2, 3\}$ and 
			${\cal A}_{12}=\{a_1\}$, 
			${\cal A}_{21}=\{a_2, a_3\}$,
			${\cal A}_{23}=\{a_ 4\}$, 
			${\cal A}_{32}=\{a_ 5\}$, 
			${\cal A}_{13}=\emptyset$, 
			${\cal A}_{31}=\{a_6\}$, 
			${\cal A}_{11}=\{a_ 7, a_8\}$, 
			say 
			$\iota_{12}(a_1)=a_2$, $\iota_{23}(a_4)=a_5$, 
			$\iota_{11}(a_7)=a_7$, $\iota_{11}(a_8)=a_8$, 
			i.e., $a_1^{-1}=a_2$, $a_4^{-1}=a_5$, $a_7^{-1}=a_7$, $a_8^{-1}=a_8$, 
			and $a_3, a_6$ have no inverse. 
			In this case, we have: 
			$
				\Phi_{\Delta}^{(1)}
				=
				\{(1, 2), (2, 3)\}
			$, 
			$
				\Phi_{\Delta}^{(2)}
				=
				\{(1, 1)\}
			$, 
			$
				\Phi_{\Delta}^{(3)}
				=
				\{(1, 3)\}
			$; 
			$
				{\cal A}^{(1)}
				=
				\{a_1, a_4\}
			$, 
			$
				{\cal A}^{(-1)}
				=
				\{a_2, a_5\}
			$, 
			$
				\overline{{\cal A}^{(1)}}
				=
				\{a_3\}
			$, 
			$
				{\cal A}^{(2)}
				=
				\{a_7, a_8\}
			$, 
			$
				{\cal A}^{(3)}
				=
				\{a_6\}
			$; 
			$c_{a_1}(t)=c_{a_2}(t)=1-\upsilon(a_1)\upsilon(a_2)t^2$, 
			$c_{a_4}(t)=c_{a_5}(t)=1-\upsilon(a_4)\upsilon(a_5)t^2$, 
			$c_{a_7}(t)=1+\upsilon(a_7)t$, 
			$c_{a_8}(t)=1+\upsilon(a_8)t$, 
			$c_{a_3}(t)=c_{a_6}(t)=1$; 
			and 
			$$
				\begin{array}{l}
				A_{\Delta}(\theta^{\rm G})
				=
				\left[
					\begin{array}{ccc}
					\frac{\tau(a_7)}{1+\upsilon(a_7)t}+\frac{\tau(a_8)}{1+\upsilon(a_8)t} 
						& \frac{\tau(a_1)}{1-\upsilon(a_1)\upsilon(a_2)t^2} & 0\\
					\frac{\tau(a_2)}{1-\upsilon(a_2)\upsilon(a_1)t^2}+\tau(a_3) & 0 
						& \frac{\tau(a_4)}{1-\upsilon(a_4)\upsilon(a_5)t^2}\\
					\tau(a_6) & \frac{\tau(a_5)}{1-\upsilon(a_5)\upsilon(a_4)t^2} & 0
					\end{array}
				\right], \\\\
				B_{\Delta}(\theta^{\rm G})
				=
				\left[
					\begin{array}{ccc}
					\frac{\tau(a_1)\upsilon(a_2)}{1-\upsilon(a_1)\upsilon(a_2)t^2} &0 & 0\\
					0 
					& \frac{\tau(a_2)\upsilon(a_1)}{1-\upsilon(a_2)\upsilon(a_1)t^2}
					+\frac{\tau(a_4)\upsilon(a_5)}{1-\upsilon(a_4)\upsilon(a_5)t^2} &0\\
					0 & 0 & \frac{\tau(a_5)\upsilon(a_4)}{1-\upsilon(a_5)\upsilon(a_4)t^2}
					\end{array}
				\right].
				\end{array}
			$$
		}
		\end{ex}

		\begin{figure}[h]
		\centering
		\includegraphics[scale=0.37]{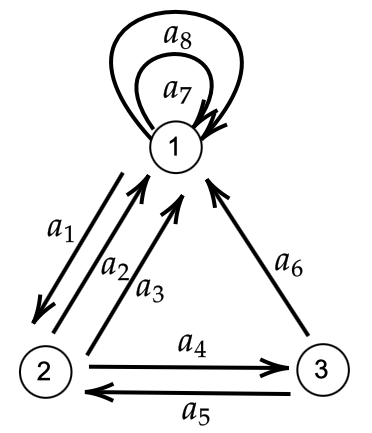} 
		\caption{The digraph $\Delta$ in Example \ref{example : Example}}
		\end{figure}


		\begin{rem}
		{\em
		Let $\Gamma=(V, E)$ be a finite simple graph and $\Delta=\Delta(\Gamma)$ 
		the symmetric digraph. 
		Note that by definition $\Gamma$ has no loops. 
		Then one can see that $A$ and $B$ are natural generalization of 
		the adjacency matrix $A_{\Gamma}$ and the degree matrix $D_{\Gamma}$ of $\Gamma$ respectively 
		(c.f., \cite{IIMSS21}). 
		In this case, we have $|{\cal A}_{uv}|= |{\cal A}_{vu}|=1$ for non-empty ${\cal A}(u, v)$. 
		In addition, if we consider the case where $\theta^{\rm G}=\theta^{\rm I}$, i.e., $\tau=\upsilon=1$, 
		then it follows that 
		$$c_a(t)=1-t^2$$ 
		for all $a\in {\cal A}$, 
		since $|{\cal A}_{uv}|=1$ if ${\cal A}_{uv}\neq\emptyset$ and 
		the inclusions $\iota_{uv} : {\cal A}_{uv}\rightarrow {\cal A}_{vu}$ are bijective. 
		A simple observation shows that, for $u, v\in V$,  
		$$
			a_{uv}
			=
			\frac{|{\cal A}_{uv}|}{1-t^2}, 
			\quad
			b_{uu}
			=
			\frac{|{\cal A}^{(1)}\cap{\cal A}_{u*}|}{1-t^2}.
		$$ 
		One can easily see that $|{\cal A}_{uv}|=1$ iff $\{u, v\}\in E$ (otherwise zero), 
		and $|{\cal A}^{(1)}\cap{\cal A}_{u*}|$ gives 
		the number of edges in $\Gamma$ satisfying $\{u, v\}\in E$ for some $v\in V$. 
		This shows that 
		$$
			A_{\Delta}(\theta^{\rm G})
			=
			\frac{1}{1-t^2}A_{\Gamma},
			\quad
			B_{\Delta}(\theta^{\rm G})
			=
			\frac{1}{1-t^2}D_{\Gamma}.
		$$
		Thus we can regard that the weighted adjacency matrix and the weighted backtrack matrix are, respectively, 
		natural generalization of the adjacency matrix and the degree matrix. 
		}
		\end{rem}

		\vspace*{5mm}
		Let $\Delta=(V, {\cal A})$ be a finite digraph. 
		Recall that $M=M_{\Delta}(\theta^{\rm G})=(\theta^{\rm G}(a, a'))_{a, a'\in{\cal A}}$. 
		Let 
		\begin{eqnarray*}
		&&H=H_{\Delta}(\theta^{\rm G})=(\tau(a')\delta_{\head(a)\tail(a')})_{a, a'\in{\cal A}},\\
		&&J=J_{\Delta}(\theta^{\rm G})=(\upsilon(a')\delta_{a^{-1}a'})_{a, a'\in{\cal A}},\\
		&&K=K_{\Delta}(\theta^{\rm G})=(\delta_{\head(a)v})_{a\in{\cal A}, v\in V},\\
		&&L=L_{\Delta}(\theta^{\rm G})=(\tau(a')\delta_{u\tail(a')})_{u\in V, a'\in{\cal A}}.
		\end{eqnarray*}
		For each arc $a\in{\cal A}$, 
		we consider the following restrictions 
		$$
		\begin{array}{l}
			J(a)=(\upsilon(a')\delta_{\alpha^{-1}\alpha'})_{\alpha, \alpha'\in{\cal E}(a)},\\ 
			K(a)=(\delta_{\head(\alpha), v})_{\alpha\in{\cal E}(a), v\in V},\\ 
			L(a)=(\tau(a')\delta_{u\tail(\alpha')})_{u\in V, \alpha'\in{\cal E}(a)}
		\end{array}
		$$
		for the matrices $J$, $K$, and $L$. 
		Note that $J(a)$ is $2\times 2$-matrix if $a\in{\cal A}^{(1)}$, 
		$1\times1$ otherwise. 
		Hence we can arrange the arcs so as to 
		the matrix $J$ is a direct sum
		$$
			J=
			\left(\bigoplus_{a\in{\cal A}^{(1)}}J(a)\right)
			\oplus 
			\left(\bigoplus_{a\in{\cal A}\setminus ({\cal A}^{(1)}\cup{\cal A}^{(-1)})}J(a)\right)
		$$ 
		of $2\times 2$ blocks and $1\times 1$ blocks. 
		We fix such a total order on ${\cal A}$. 
		If we denote by $I(a)$ ($a\in{\cal A}$) the identity matrix of degree $|{\cal E}(a)|$, 
		then the matrix $I+tJ$ is a direct sum of the matrices 
		$\oplus_{a\in{\cal A}\setminus{\cal A}^{(-1)}}(I(a)+tJ(a))$, 
		where the direct summands are all invertible on $R[[t]]$.

		\begin{lem}
		The matrix $I+tJ$ is invertible.
		\end{lem}

		For $\Delta$ and $\theta^{\rm G}$, 
		we denote by $I_{\Delta}(t; \theta^{\rm G})$ the following formal power series with indeterminate $t$:
		$$
			\frac
			{1}
			{\det(I+tJ)\det(I-tA_{\Delta}(\theta^{\rm G})+t^2B_{\Delta}(\theta^{\rm G}))}.
		$$
		\begin{thm}[Main theorem]
		Let $\Delta$, $R$ and $\theta^{\rm G}$ be as above. 
		We have $$Z_{\Delta}(t; \theta^{\rm G})=I_{\Delta}(t; \theta^{\rm G}).$$
		\end{thm}
		{\it Proof.}
		Recall that we have the identity
		$$
			Z_{\Delta}(t; \theta^{\rm G})
			=
			\frac{1}{\det(I-tM)}, 
		$$
		where $M=M_{\Delta}(\theta^{\rm G})$. 
		Let $H, J, K$ and $L$ be as above. 
		By definition, it follows that $M=H-J$. 
		It also follows that $H=KL$, thus $M=KL-J$. 
		Hence we have 
		\begin{eqnarray*}
		\det(I-tM)
		&=&
		\det(I-t(KL-J))\\
		&=&
		\det((I+tJ)-tKL)\\
		&=&
		\det(I+tJ)\det(I-t(I+tJ)^{-1}KL)\\
		&=&
		\det(I+tJ)\det(I-tL(I+tJ)^{-1}K),
		\end{eqnarray*}
		where the final identity follows from 
		the well-known identity $\det(I-AB)=\det(I-BA)$ in linear algebra. 
		Since each direct summand of 
		$$
			I+tJ
			=
			\bigoplus_{a\in {\cal A}\setminus{\cal A}^{(-1)}}
			I(a)+tJ(a)
		$$ 
		is invertible, we have 
		$
			(I+tJ)^{-1}
			=
			\bigoplus_{a\in {\cal A}\setminus{\cal A}^{(-1)}}
			(I(a)+tJ(a))^{-1},
		$ 
		and it follows that 
		$$
			L(I+tJ)^{-1}K
				=
				\sum_{a\in{\cal A}\setminus{\cal A}^{(-1)}}
				L(a)(I(a)+tJ(a))^{-1}K(a).
		$$
		Note that $\det(I(a)+tJ(a))=c_a(t)$ for $a\in{\cal A}\setminus{\cal A}^{(-1)}$, 
		and we have 
		$$
			(I(a)+tJ(a))^{-1}
			=
			\left\{
				\begin{array}{ll}
				c_a(t)^{-1}(I(a)-tJ(a)), & \mbox{if $a\in{\cal A}^{(1)}$}\\
				c_a(t))^{-1}I(a), & \mbox{if $a\in{\cal A}^{(2)}\sqcup{\cal A}^{\times}$}.
				\end{array}
			\right.
		$$
		Hence it follows that 
		$$
			\sum_{a\in{\cal A}\setminus{\cal A}^{(-1)}}
			L(a)(I(a)+tJ(a))^{-1}K(a)
			=
			\sum_{a\in{\cal A}\setminus{\cal A}^{(-1)}}
			c_a(t)^{-1}L(a)K(a)
			-t
			\sum_{a\in{\cal A}^{(1)}}
			c_a(t)^{-1}L(a)J(a)K(a). 
		$$

		The $(u, v)$-entry $r_{uv}$ of the matrix 
		$
		\sum_{a\in{\cal A}\setminus{\cal A}^{(-1)}}
		c_a(t)^{-1}L(a)K(a)
		$ 
		is given by 
		\begin{equation}\label{(u, v)-entry}
		r_{uv}=
		\sum_{a\in{\cal A}\setminus{\cal A}^{(-1)}}
		c_a(t)^{-1}
		\sum_{\alpha\in{\cal E}(a)}
		\tau(\alpha)\delta_{u\tail(\alpha)}\delta_{\head(\alpha)v}.
		\end{equation}
		Note that 
		$
		\delta_{u\tail(\alpha)}\delta_{\head(\alpha)v}\neq 0
		$ is equivalent to 
		$a\in{\cal A}_{uv}$. 
		It follows that 
		$$
			r_{uv}
			=
			\sum_{a\in{\cal A}\setminus{\cal A}^{(-1)}}
			c_a(t)^{-1}
			\sum_{\alpha\in{\cal E}(a)\cap{\cal A}_{uv}}
			\tau(\alpha).
		$$
		We verify $r_{uv}=a_{uv}$ for all $(u, v)\in V\times V$. 
		Suppose that $u\preceq v$. 
		If $a\in {\cal A}^{{(1)}}$, 
		then 
		${\cal E}(a)\cap{\cal A}_{uv}=\{a\}$. 
		Otherwise, 
		we have 
		$
			{\cal E}(a)\cap{\cal A}_{uv}=\emptyset.
		$ 
		This implies that 
		$
			r_{uv}
			=
			\sum_{a\in{\cal A}\setminus{\cal A}^{(-1)}}
			c_a(t)^{-1}\tau(a).
		$ 
		In the case where $v\prec u$, 
		${\cal E}(a)\cap{\cal A}_{uv}\neq\emptyset$ 
		implies that 
		$a\in{\cal A}^{\times}$, 
		and we have ${\cal E}(a)\cap{\cal A}_{uv}=\{a\}$ 
		for $a\in{\cal A}^{\times}$. 
		Suppose that $u=v$. 
		In this case, 
		${\cal E}(a)\cap{\cal A}_{uu}\neq\emptyset$ 
		implies $a\in{\cal A}^{(2)}$, 
		and we have ${\cal E}(a)\cap{\cal A}_{uu}=\{a\}$. 
		Therefore, putting all these together, 
		it follows that $r_{uv}=a_{uv}$ for all $(u, v)\in V\times V$.

		The $(u, v)$-entry $s_{uv}$ of the matrix 
		$
		\sum_{a\in{\cal A}^{(1)}}
		c_a(t)^{-1}L(a)J(a)K(a)
		$ 
		is given by 
		\begin{equation}\label{s_{uv}}
		s_{uv}=
		\sum_{a\in{\cal A}^{(1)}}
		c_a(t)^{-1}
		\sum_{\alpha, \beta\in{\cal E}(a)}
		\tau(\alpha)\upsilon(\beta)
		\delta_{u\tail(\alpha)}
		\delta_{\alpha^{-1}\beta}
		\delta_{\head(\beta)v}.
		\end{equation}
		We verify $s_{uv}=b_{uv}$ for any $(u, v)\in V\times V$. 
		Let $a\in{\cal A}^{(1)}$. 
		We have ${\cal E}(a)=\{a, a^{-1}\}$ with $a\neq a^{-1}$. 
		Thus it follows that  
		$$
			s_{uv}
			=
			\sum_{a\in{\cal A}^{(1)}}
			c_a(t)^{-1}
			\tau(a)\upsilon(a^{-1})
			\delta_{u\tail(a)}
			\delta_{\head(a^{-1})u},
		$$
		which equals 
		$
			\sum_{a\in{\cal A}^{(1)}\cap{\cal A}_{u*}}
			c_a(t)^{-1}
			\tau(a)\upsilon(a^{-1}).
		$
		Now we have show that $s_{uv}=b_{uv}$. \qed

		\vspace*{5mm}
		\noindent
		{\bf Acknowledgements}
		
		The authors would like to express his deep gratitude to Professor Iwao Sato, 
		Oyama National College of Technology, 
		who suggested the problem considered in this article, 
		for illuminating discussions and valuable comments. 
		The authors also would like to thank anonymous referees 
		for his valuable comments that improve the article. 
		The first named author is partially supported by Grant-in-Aid for JSPS Fellows, 
		Grant Number JP20J20590. 
		The second named author is partially supported by JSPS KAKENHI, 
		Grant Number JP22K03262. 
		The authors also would like to thank anonymous referees 
		for his valuable comments that improve the article.

{

\end{document}

		\begin{rem}
		{\em 
		Let $\Delta$ be the symmetric digraph of a finite graph. 
		If $\tau=\upsilon=1$, 
		then $Z_{\Delta}(t; \theta^{\rm G})$ is nothing but the Ihara zeta function
		\cite{bass92, hashimoto90, ihara66, kotanisunada00, serre80}. 
		If $\tau=\upsilon$, it is the graph zeta treated in \cite{mizunosato04}, 
		which is called the {\it Mizuno-Sato zeta function}. 
		If $\upsilon=1$, it is the one defined in \cite{sato07}, 
		called the {\it Sato zeta function}. 
		Let $q$ be an indeterminate, and replace $R$ by the polynomial ring $R[q]$. 
		If we let $\tau=1$ and $\upsilon=1-q$, 
		then the resulting graph zeta is called the {\it Bartholdi zeta function} defined in \cite{bartholdi99}. 
		The {\it egde zeta function} and the {\it path zeta function}\cite{starkterras96} also come out from this framework 
		\cite{morita20}.  
		}
		\end{rem}


		\subsection{Dynamical systems on finite digraphs}
		
		Let $\Delta=(V, {\cal A})$ be a finite digraph. 
		Consider the set $\Xi$ of two-sided infinite sequences 
		$(a_i)_{i\in{\mathbb Z}}\in {\cal A}^{\mathbb Z}$ on ${\cal A}$ 
		satisfying $\head(a_i)=\tail (a_{i+1})$ for all $i\in{\mathbb Z}$, 
		and the left shift operator 
		$
			\varphi : {\cal A}^{\mathbb Z}\rightarrow {\cal A}^{\mathbb Z} 
			: (a_i)_i\mapsto (a_{i+1})_i
		$. 
		If $x\in X_m$, then the integer $m$ is called a {\it period} of $x$. 
		Note that a period of $x\in X$ is not uniquely determined for $x$, 
		since any multiple of a period of $x$ is also a period. 
		Since the restriction $\lambda=\varphi|_{\Xi}$ gives a bijection on $\Xi$, 
		one has a dynamical system $(\Xi, \lambda)$. 
		Let $X_m={\rm Fix}(\lambda^m)$ denote the set of $m$-{\it periodic point}, 
		i.e., an element $x\in \Xi$ satisfying $\lambda^m(x)=x$, 
		and $X=\cup_{m\geq 1}X_m$ the set of periodic point of $(\Xi, \lambda)$. 
		An element of $X_m$ is also called a {\it closed path} of {\it length} $m$ in $\Delta$. 
		Let $x\in X$. 
		The minimum period 
		$$
			\varpi(x)
			=
			\min\{m\mid x\in X_m\}
		$$ 
		is called the {\it prime period} of $x$. 
		It is obvious that $x\in X_{\varpi(x)}$ for any $x\in X$. 
		If we regard $x\in X$ as an element of $X_{\varpi(x)}$, 
		then $x\in X_{\varpi(x)}$ is called the {\it prime section} of $x$, 
		denoted by $\pi(x)$. 
		A straightforward argument shows the following lemma.

		\begin{lem}\label{lemma : the prime period divides a period}
		If $x\in X_m$, then the prime period $\varpi(x)$ divides $m$. 
		\end{lem}

{}